\documentclass[12pt,a4paper]{article}
\usepackage{amsmath,amsfonts,amssymb,amscd}

\def\gr{\operatorname{gr}}

\def\id{\operatorname{id}}

\def\GL{\operatorname{GL}}

\def\Aut{\operatorname{Aut}}

\newcounter{th}
\def\t{\refstepcounter{th}{\bf \noindent{Theorem} \arabic{th}. }}

\newcounter{le}
\def\l{\refstepcounter{le}{\bf \noindent{Proposition} \arabic{le}. }}

\newcounter{lem}

\newcounter{de}

\newcounter{ex}

\begin{document}

\begin{center}
{\Large Even-homogeneous supermanifolds on the complex projective
line} \footnote[1]{Supported by the Mathematisches
Forschungsinstitut Oberwolfach, Leibniz Fellowship}

\medskip

{\bf E.G.~Vishnyakova}\\
\end{center}

\begin{quotation}
\small \noindent\textsc{Abstract.} The classification of
even-homogeneous complex supermanifolds of dimension $1|m$, $m\leq
3$, on $\mathbb{CP}^1$ up to isomorphism is given. An explicit
description of such supermanifolds in terms of local charts and
coordinates is obtained.

\end{quotation}

\bigskip
{\bf 1. Introduction.} The study of homo\-ge\-neous supermanifolds
with under\-lying manifold $\mathbb {CP}^1$ was started in
\cite{Bunegina_Oni}. There the classification of homogeneous complex
supermanifolds of dimension $1|m$, $m\le 3$, up to isomorphism was
given. The purpose of this paper is to classify up to isomorphism
even-homo\-ge\-neous non-split complex supermanifolds of dimension
$1|m$, $m\le 3$, on $\mathbb {CP}^1$. Some other classification
results concerning non-split complex supermanifolds on $\mathbb
{CP}^n$ can be found in \cite{Bash,OniP1,OniP2}. 

The paper is structured as follows. In Section $2$ we explain the
idea of the classification. Similar idea was used in
\cite{Bunegina_Oni} by the classification of homogeneous
supermanifolds on $\mathbb{CP}^1$. In Section $3$ we calculate the
$1$-cohomology group with values in the tangent sheaf. We use here
an easier way than in \cite{Bunegina_Oni}, which permits to classify
even-homogeneous supermanifolds.

 By the Green Theorem
we can assign a supermanifold to each cohomology class of the
$1$-cohomology group. In Section $4$ we find out cohomo\-logy
classes corresponding to even-homogeneous supermanifolds. Notice
that these supermanifolds can be isomorphic. The classification up
to isomorphism of even-homogeneous complex supermanifolds of
dimension $1|m$, $m\le 3$, on $\mathbb{CP}^1$ is obtained in Section
$5$.

\bigskip
{\bf 2. Even-homogeneous supermanifolds on $\mathbb {CP}^1$.} We
study complex analytic supermanifolds in the sense of
\cite{Bunegina_Oni,Man}. If $\mathcal{M} =(\mathcal{M}_0,{\mathcal
O}_{\mathcal{M}})$ is a supermanifold, we denote by $\mathcal{M}_0$
the {\it underlying complex manifold} of $\mathcal{M}$ and by
${\mathcal O}_{\mathcal{M}}$ the {\it structure sheaf} of
$\mathcal{M}$, i.e. the sheaf of commutative associative complex
superalgebras on $\mathcal{M}_0$. Denote by ${\mathcal
T}_{\mathcal{M}}$ the {\it tangent sheaf} of $\mathcal{M}$, i.e. the
sheaf of derivations of the structure sheaf
$\mathcal{O}_{\mathcal{M}}$. Denote by $({\mathcal
T}_{\mathcal{M}})_{\bar 0}\subset {\mathcal T}_{\mathcal{M}}$ the
subsheaf of all even vector fields. An {\it action of a Lie group
$G$ on a supermanifold} $\mathcal{M}$ is a morphism $\nu =
(\nu_0,\nu^*):G \times \mathcal{M}\to \mathcal{M}$ such that it
satisfies the usual conditions, modeling the action axioms. An
action $\nu$ is called {\it even-transitive} if $\nu_{0}$ is
transitive.
 A supermanifold $\mathcal{M}$
is called {\it even-homogeneous} if it possesses an even-transitive
action of a Lie group.

Assume that $\mathcal{M}_0$ is compact and connected. It is
well-known that the group of all automorphisms of $\mathcal{M}$,
which we denote by $\Aut \mathcal{M}$, is a Lie group with the Lie
algebra $H^0(\mathcal{M}_0, (\mathcal{T}_{\mathcal{M}})_{\bar 0})$.
(Recall that by definition any morphism of a super\-mani\-fold is
even.) Let us take any homomorphism of Lie algebras $\varphi:
\mathfrak{g} \to H^0(\mathcal{M}_0,(\mathcal{T}_{\mathcal{M}})_{\bar
0})$. We can assign the homomorphism of Lie groups $\Phi:G\to \Aut
\mathcal{M}$ to $\varphi$, where $G$ is the simple connected Lie
group with the Lie algebra $\mathfrak{g}$. Notice that $\Phi$ is
even-transitive iff the image of $\mathfrak{g}$ in
$H^0(\mathcal{M}_0, (\mathcal{T}_{\mathcal{M}})_{\bar 0})$ generates
the tangent space $T_x(\mathcal{M})$ at any point $x\in
\mathcal{M}_0$.

In this paper we will consider the case $\mathcal{M}_0=\mathbb
{CP}^1$. Therefore, the classification problem reduces to the
following problem: \textit{to classify up to isomorphism complex
supermanifolds $\mathcal{M}$ of dimension $1|m$, $m\le 3$, such that
$H^0(\mathcal{M}_0,(\mathcal{T}_{\mathcal{M}})_{\bar 0})$ generates
the tangent space $T_x(\mathcal{M})$ at any point $x\in
\mathcal{M}_0$}.

 Recall that a supermanifold $\mathcal{M}$ is
called {\it split} if ${\mathcal
O}_{\mathcal{M}}\simeq\bigwedge{\mathcal E}$, where ${\mathcal E}$
is a sheaf of sections of a vector bundle ${\mathbf E}$ over
$\mathcal{M}_0$. In this case $\dim\,\mathcal{M} = n|m$, where $n =
\dim \mathcal{M}_0$ and $m$ is the rank of ${\mathbf E}$. The
structure sheaf ${\mathcal O}_{\mathcal{M}}$ of a split
supermanifold possesses by definition the ${\mathbb Z}$-grading; it
induces the ${\mathbb Z}$-grading in ${\mathcal T}_{\mathcal{M}} =
\bigoplus_{p=-1}^m({\mathcal T}_{\mathcal{M}})_p$. Hence, the
superspace $H^0(\mathcal{M}_0,{\mathcal T}_{\mathcal{M}})$ is also
${\mathbb Z}$-graded. Consider the subspace ${\mathrm
{End}}\,{\mathbf E}\subset H^0(\mathcal{M}_0,{\mathcal
T}_{\mathcal{M}})_0$ consisting of all endo\-mor\-phisms of the
vector bundle ${\mathbf E}$, which induce the identity morphism on
$\mathcal{M}_0$. Denote by ${\mathrm {Aut}}\,{\mathbf
E}\subset{\mathrm {End}}\,{\mathbf E}$ the group of automorphisms
containing in ${\mathrm {End}}\,{\mathbf E}$. We define an action
${\mathrm{Int}}$ of ${\mathrm {Aut}}\,{\mathbf E}$ on ${\mathcal
T}_{\mathcal{M}}$ by ${\mathrm{Int}} A: v\mapsto A v A^{-1}$. Since
the action preserves the $\mathbb{Z}$-grading, we have the action of
${\mathrm {Aut}}\,{\mathbf E}$ on $H^1(\mathcal{M}_0,{(\mathcal
T}_{\mathcal{M}})_2)$.

We can assign the split supermanifold $\gr \mathcal{M}=
(\mathcal{M}_0, \mathcal{O}_{\gr\mathcal{M}})$ to each supermanifold
$\mathcal{M}$, see e.g. \cite{Bunegina_Oni}. It is called the {\it
retract} of $\mathcal{M}$. To classify supermanifolds, we will use
the following corollary of the well-known Green Theorem (see e.g.
\cite{Bunegina_Oni} for more details).

\medskip

\t\label{Green corollary}[Green] {\it Let $\widetilde{\mathcal{M}} =
(\mathcal{M}_{0}, \bigwedge \mathcal{E})$ be a split supermanfold of
dimension $n|m$, where $m\le 3$. Then classes of isomorphic
supermanifolds $\mathcal{M}$ with the retract $\gr \mathcal{M} =
\widetilde{\mathcal{M}}$ are in bijection with orbits of the action
${\mathrm{Int}}$ of the group ${\mathrm {Aut}}\,{\mathbf E}$ on
$H^1(\mathcal{M}_0,{(\mathcal T}_{\widetilde{\mathcal{M}}})_2)$.}

\medskip

\noindent {\bf Remark.} This theorem permits to classify
supermanifolds $\mathcal{M}$ such that $\gr \mathcal{M}$ is fix up
to isomorphisms which induce identity morphism on $\gr \mathcal{M}$.

\medskip

In what follows we will consider the case
$\mathcal{M}_0={\mathbb{CP}}^1$. Let $\mathcal{M}$ be a
super\-mani\-fold of dimension $1|m$. Denote by $U_0$ and $U_1$ the
standard charts on ${\mathbb{CP}}^1$ with coordinates  $x$ and
$y=\frac1x$ respectively. By the Grothendieck Theorem we can cover
$\gr\mathcal{M}$ by two charts $(U_0,\mathcal{O}_{\gr
\mathcal{M}}|_{U_0})$ and $(U_1,\mathcal{O}_{\gr
\mathcal{M}}|_{U_1})$ with local coordinates
$x,\,\xi_1,\ldots,\xi_m$ and $y,\,\eta_1,\ldots,\eta_m$,
respectively, such that in $U_0\cap U_1$ we have
$$
y = x^{-1},\quad \eta_i = x^{-k_i}\xi_i,\; i = 1,\ldots,m,
$$
where $k_i,\; i = 1,\ldots,m$, are integers.  We will identify
$\gr\mathcal{M}$ with the set $(k_1,\ldots, k_m)$. Note that a
permutation of $k_i$ induces the automorphism of $\gr\mathcal{M}$.
It was shown that any supermanifold $\gr\mathcal{M}$ is
even-homogeneous, see \cite{Bunegina_Oni}, Formula $(18)$. The
following theorem was also proved in \cite{Bunegina_Oni},
Proposition $14$:

\medskip

\t\label{even-homogen kriterij}{\it Assume that $m \leq 3$ and
$\mathcal{M}_0= {\mathbb{CP}}^1$. Let $\mathcal{M}$ be a
supermanifold with the retract $\gr\mathcal{M}=\bigwedge
\mathcal{E}$, which corresponds to the cohomology class $\gamma\in
H^1(\mathcal{M}_0,({\mathcal T}_{\gr\mathcal{M}})_2)$ by Theorem
\ref{Green corollary}. The following conditions are equivalent:

\smallskip

$1.$ The supermanifold $\mathcal{M}$ is even-homogeneous.

$2.$ There is a subalgebra ${\mathfrak a}\simeq {\mathfrak
{sl}}_2({\mathbb C})$ such that
\begin{equation}\label{subalgebras a}
H^0(\mathcal{M}_0,({\mathcal T}_{\gr\mathcal{M}})_0)=
{\mathrm{End}}\,{\mathbf E}\supset\!\!\!\!\!\!\!+{\mathfrak a},
\end{equation}
and $[v,\gamma]=0$ in $H^1(\mathcal{M}_0,({\mathcal
T}_{\gr\mathcal{M}})_2)$ for all $v\in{\mathfrak a}$}.

\medskip

Here ${\mathbf E}$ is the vector bundle corresponding to the locally
free sheaf $\mathcal{E}$.

From now on we will omit the index $\gr\mathcal{M}$ and will denote
by ${\mathcal T}$ the sheaf of derivations of $\mathcal{O}_{\gr
\mathcal{M}}$. Recall that the sheaf $\mathcal{O}_{\gr \mathcal{M}}$
is $\mathbb{Z}$-graded; it induces the $\mathbb{Z}$-grading in
$\mathcal{T} = \bigoplus_p\mathcal{T}_p$. Denote by
$H^1({\mathbb{CP}}^1,{\mathcal T}_2)^{\mathfrak a}\subset
H^1({\mathbb{CP}}^1,{\mathcal T}_2)$ the subset of ${\mathfrak
a}$-invariants, i.e. the set of all elements $w$ such that $[v,w]=0$
for all $v\in{\mathfrak a}$. The supermanifold corresponding to a
cohomology class $\gamma\in H^1({\mathbb{CP}}^1,{\mathcal
T}_2)^{\mathfrak a}$ by Theorem \ref{Green corollary} is called {\it
${\mathfrak a}$-even-homogeneous.}

The description of subalgebras ${\mathfrak a}$ satisfying
(\ref{subalgebras a}) up to conjugation by elements from
${\mathrm{Aut}}\,{\mathbf E}$ and up to renumbering of $k_i$ was
obtained in \cite{Bunegina_Oni}:

\smallskip
\noindent
1) ${\mathfrak a}={\mathfrak s}=
\langle {\mathbf e}=\frac{\partial}{\partial x},
{\mathbf f}=\frac{\partial}{\partial y},{\mathbf h}=
[{\mathbf e},{\mathbf f}] \rangle$.

\noindent 2) ${\mathfrak a}={\mathfrak s'}=\langle {\mathbf
e'}=\frac{\partial}{\partial x}+\xi_2\frac{\partial}{\partial
\xi_1}, {\mathbf f'}=\frac{\partial}{\partial y}+
\eta_1\frac{\partial}{\partial \eta_2}, {\mathbf h'}=[{\mathbf
e'},{\mathbf f'}]\rangle$ if $k_1=k_2$.

\noindent 3) ${\mathfrak a}={\mathfrak s''}=\langle {\mathbf
e''}=\frac{\partial}{\partial x}+ \xi_2\frac{\partial}{\partial
\xi_1}+\xi_3\frac{\partial}{\partial \xi_2}, {\mathbf
f''}=\frac{\partial}{\partial y}+ 2\eta_1\frac{\partial}{\partial
\eta_2}+ 2\eta_2\frac{\partial}{\partial \eta_3}, {\mathbf
h''}=[{\mathbf e''},{\mathbf f''}]\rangle$ if $k_1=k_2=k_3$.

\bigskip

{\bf 3. Basis of $H^1({\mathbb{CP}}^1,{\mathcal T}_2).$} Assume that
$m = 3$. Let $\mathcal{M}$ be a split super\-mani\-fold,
$\mathcal{M}_0=\mathbb{CP}^1$ be its reduction and ${\mathcal T}$ be
its tangent sheaf. In \cite{Bunegina_Oni} the ${\mathfrak
s}$-invariant decomposition
\begin{equation}\label{s-invar decomp}
{\mathcal T}_2=\sum_{i<j}{\mathcal T}_2^{ij}
\end{equation}
was obtained. The sheaf ${\mathcal T}_2^{ij}$ is a locally free
sheaf of rank $2$; its basis sections over
$(U_0,\mathcal{O}_{\mathcal{M}}|_{U_0})$ are:
\begin{equation}\label{loc section of T^ij}
\xi_i\xi_j\frac{\partial}{\partial x}, \ \
\xi_i\xi_j\xi_l\frac{\partial}{\partial \xi_l};
\end{equation}
where $l\ne i,j$. In $U_0\cap U_1$ we have
\begin{equation}\label{sootnoshen v peresech}
\begin{aligned}
\xi_i\xi_j\frac{\partial}{\partial x}
&=-y^{2-k_i-k_j}\eta_i\eta_j\frac{\partial}{\partial y} -
k_l y^{1-k_i-k_j}\eta_i\eta_j\eta_l\frac{\partial}{\partial \eta_l},\\
\xi_i\xi_j\xi_l\frac{\partial}{\partial \xi_l}
&= y^{-k_i-k_j}\eta_i\eta_j\eta_l\frac{\partial}{\partial \eta_l}.
\end{aligned}
\end{equation}

Let us calculate a basis of $H^1({\mathbb{CP}}^1,{\mathcal
T}_2^{ij})$. We will use the \v{C}ech cochain complex of the cover
${\mathfrak U} =\{U_0,U_1\}$. Hence, $1$-cocycle with values in the
sheaf ${\mathcal T}_2^{ij}$ is a section $v$ of ${\mathcal
T}_2^{ij}$ over $U_0\cap U_1$. We are looking for {\it basis
cocycles}, i.e. cocycles such that their cohomology classes form a
basis of
 $H^1({\mathfrak U},{\mathcal
T}_2^{ij})\simeq H^1({\mathbb{CP}}^1,{\mathcal T}_2^{ij})$. Note
that if $v\in Z^1({\mathfrak U},{\mathcal T}_2^{ij})$ is holomorphic
in $U_0$ or $U_1$ then the cohomology class of $v$ is equal to $0$.
Obviously, any $v\in Z^1({\mathfrak U},{\mathcal T}_2^{ij})$ is a
linear combination of vector fields (\ref{loc section of T^ij}) with
holomorphic in $U_0\cap U_1$ coefficients. Further, we expand these
coefficients in a Laurent series in $x$ and drop the summands
$x^n,\; n\ge 0$, because they are holomorphic in $U_0$. We see that
$v$ can be replaced by
\begin{equation}\label{cocycle chomological to v}
v = \sum_{n=1}^{\infty} a^n_{ij}x^{-n}\xi_i\xi_j\frac{\partial}
{\partial x} + \sum_{n=1}^{\infty}
b^n_{ij}x^{-n}\xi_i\xi_j\xi_l\frac{\partial} {\partial \xi_l},
\end{equation}
where $a^n_{ij}, b^n_{ij}\in{\mathbb C}$. Using (\ref{sootnoshen v
peresech}), we see that the summands corresponding to $n\ge
k_i+k_j-1$ in the first sum of (\ref{cocycle chomological to v}) and
the summands corresponding to $n\ge k_i+k_j$ in the second sum of
(\ref{cocycle chomological to v}) are holomorphic in $U_1$. Further,
it follows from (\ref{sootnoshen v peresech}) that
$$
x^{2-k_i-k_j}\xi_i\xi_j\frac{\partial}{\partial x}\sim
-k_lx^{1-k_i-k_j}\xi_i\xi_j\xi_l\frac{\partial}{\partial \xi_l}.
$$
Hence the cohomology classes of the following cocycles
\begin{equation}\label{cocycles generating H^1}
\begin{aligned}
&x^{-n}\xi_i\xi_j\frac{\partial}{\partial x},\ \ n=1,{\ldots},k_i+k_j-3,\\
&x^{-n}\xi_i\xi_j\xi_l\frac{\partial}{\partial \xi_l},\ \
n=1,{\ldots},k_i+k_j-1,
\end{aligned}
\end{equation}
generate $H^1({\mathbb{CP}}^1,{\mathcal T}_2^{ij})$. If we examine
linear combination of (\ref{cocycles generating H^1}) which are
cohomological trivial, we get the following theorem.

\medskip

\t\label{basis of H^1}{\it Assume that $i<j$, $l\ne i,j$. The basis
of $H^1({\mathbb{CP}}^1,{\mathcal T}_2^{ij})$

$1.$ is given by $(\ref{cocycles generating H^1})$ if $k_i+k_j>3$;

$2.$ is given by
$$
x^{-1}\xi_i\xi_j\xi_l\frac{\partial}{\partial \xi_l},\quad
x^{-2}\xi_i\xi_j\xi_l\frac{\partial}{\partial \xi_l},
$$
if $k_i+k_j=3$;

$3.$ is given by
$$
x^{-1}\xi_i\xi_j\xi_l\frac{\partial}{\partial \xi_l},
$$
if $k_i+k_j=2$, $k_l=0$.

$4.$ If $k_i+k_j=2$, $k_l\ne 0$ or $k_i+k_j<2$, we have
$H^1({\mathbb{CP}}^1,{\mathcal T}_2^{ij})=\{0\}$.}

\medskip

Note that the similar method can be used for computation of a basis
of $H^1({\mathbb{CP}}^1,{\mathcal T}_q)$ for any $m$ and $q$.

\smallskip
{\bf 4. Basis of $H^1({\mathbb{CP}}^1,{\mathcal T}_2)^{\mathfrak
a}$}. Let us calculate a basis of $H^1({\mathbb{CP}}^1,{\mathcal
T}_2)^{\mathfrak s}$. The decomposition (\ref{s-invar decomp}) is
${\mathfrak s}$-invariant, hence,
$$
H^1({\mathbb{CP}}^1,{\mathcal T}_2)^{\mathfrak s}=
\bigoplus_{i<j} H^1({\mathbb{CP}}^1,{\mathcal T}_2^{ij})^{\mathfrak s}.
$$
Denote by $[z]$ the cohomology class corresponding to a $1$-cocycle
$z$.

\medskip

\t\label{basis of H^1(T)^a}{\it Let us fix $i<j$ and $l\ne i,j$.
Then

$1)\; H^1({\mathbb{CP}}^1,{\mathcal T}_2^{ij})^{\mathfrak s}=
\langle [\frac1x\xi_i\xi_j\frac{\partial}{\partial x}+
\frac{k_l}{2x^2}\xi_i\xi_j\xi_l\frac{\partial}{\partial
\xi_l}]\rangle$ if $k_i+k_j=4$,

$2)\; H^1( {\mathbb{CP}}^1,{\mathcal T}_2^{ij})^{\mathfrak s}=
\langle [\frac{1}{x}\xi_i\xi_j\xi_l\frac{\partial}{\partial
\xi_l}]\rangle$ if $k_i+k_j=2,\,k_l=0$,

$3)\;H^1({\mathbb{CP}}^1,{\mathcal T}_2^{ij})^{\mathfrak s}=\{0\}$
otherwise}.

\medskip

\noindent{\it Proof.} We have to find out highest vectors of the
${\mathfrak s}$-module $H^1({\mathbb{CP}}^1,{\mathcal T}_2^{ij})$
having weight $0$. By Propositions $8$ and $9$ of
\cite{Bunegina_Oni}, any cocycle $z$ from the Theorem \ref{basis of
H^1} fulfils the condition $[{\mathbf h},z] = \lambda z$. More
precisely, $\lambda = 0$ if $z =
x^{-r}\xi_i\xi_j\frac{\partial}{\partial x},\; 2r=k_i+k_j-2$ и $z =
x^{-r}\xi_i\xi_j\xi_l\frac{\partial}{\partial \xi_l},\; 2r=k_i+k_j$.
If we examine linear combination $w$ of these cocycles such that
$[{\mathbf e},w]\sim 0$, we obtain the result of the Theorem.$\Box$

\medskip

\t\label{basis of H^1(T)^s'}{\it Assume that
$H^1({\mathbb{CP}}^1,{\mathcal T}_2)^{\mathfrak s'}\ne 0$. Then we
have the following possibilities:

\noindent $1)\;(k_1,k_2,k_3)=(2,2,1)$ and a basis of
$H^1({\mathbb{CP}}^1,{\mathcal T}_2)^{\mathfrak s'}$ is given by
\begin{equation}\label{basis s'1}
[\frac{1}{x}\xi_1\xi_2\frac{\partial}{\partial x}+
\frac{1}{2x^2}\xi_1\xi_2\xi_3\frac{\partial}{\partial \xi_3}],\;\;
[\frac{1}{x^2}\xi_1\xi_2\xi_3\frac{\partial}{\partial \xi_2}-
\frac{1}{x}\xi_1\xi_2\xi_3\frac{\partial}{\partial \xi_1}];
\end{equation}

\noindent $2)\;(k_1,k_2,k_3)=(2,2,3)$ and a basis of
$H^1({\mathbb{CP}}^1,{\mathcal T}_2)^{\mathfrak s'}$ is given by
\begin{equation}\label{basis s'2}
\begin{array}{c}
[\frac{1}{x}\xi_1\xi_2\frac{\partial}{\partial x}+
\frac{3}{2x^2}\xi_1\xi_2\xi_3\frac{\partial}{\partial \xi_3}],\\

[\frac{1}{x}\xi_2\xi_3\frac{\partial}{\partial x}+
\frac{1}{x^2}\xi_1\xi_3\frac{\partial}{\partial x}+
\frac{2}{3x^2}\xi_1\xi_2\xi_3\frac{\partial}{\partial \xi_1}
-\frac{4}{3x^3}\xi_1\xi_2\xi_3\frac{\partial}{\partial \xi_2}];\\

\end{array}
\end{equation}

\noindent $3)\;(k_1,k_2,k_3)=(2,2,k_3)$, $k_3\ne 1,3$;
$(k_1,k_2,k_3)=(k,k,3-k)$, $k\ne 2$ or $(k_1,k_2,k_3)=(k,k,5-k)$,
$k\ne 2$ or $(k_1,k_2,k_3)=(1,1,0)$. Then
$$
\dim
H^1({\mathbb{CP}}^1,{\mathcal T}_2)^{\mathfrak s'}=1
$$ and a basis of
$H^1({\mathbb{CP}}^1,{\mathcal T}_2)^{\mathfrak s'}$ is given by the
following cocycles:
\begin{equation}\label{basis s'3}
\begin{array}{c}
[\frac{1}{x}\xi_1\xi_2\frac{\partial}{\partial x}+
\frac{k_3}{2x^2}\xi_1\xi_2\xi_3\frac{\partial}{\partial \xi_3}],\,\,
[\frac{1}{x^2}\xi_1\xi_2\xi_3\frac{\partial}{\partial
\xi_2} - \frac{1}{x}\xi_1\xi_2\xi_3\frac{\partial}{\partial \xi_1}],\\

[\frac{1}{x}\xi_2\xi_3\frac{\partial}{\partial x}+
\frac{1}{x^2}\xi_1\xi_3\frac{\partial}{\partial x}+
\frac{k}{3x^2}\xi_1\xi_2\xi_3\frac{\partial}{\partial
\xi_1} -\frac{2k}{3x^3}\xi_1\xi_2\xi_3\frac{\partial}{\partial \xi_2}],\\

[\frac{1}{x}\xi_1\xi_2\xi_3\frac{\partial}{\partial \xi_3}],
\end{array}
\end{equation}
respectively.

}

\medskip

\noindent{\it Proof.} Use similar argument as in Theorem \ref{basis
of H^1(T)^a}.$\Box$

\medskip

The calculation of $H^1({\mathbb{CP}}^1,{\mathcal T}_2)^{\mathfrak
s}$ and $H^1({\mathbb{CP}}^1,{\mathcal T}_2)^{\mathfrak s'}$ was
already done in \cite{Bunegina_Oni}, Proposition $19$ and
Proposition $21$, using more difficult methods. Note that the case
$2$ of Theorem \ref{basis of H^1(T)^a} and the case
$(k_1,k_2,k_3)=(1,1,0)$ of Theorem \ref{basis of H^1(T)^s'} was lost
in \cite{Bunegina_Oni}. Furthermore, in \cite{Bunegina_Oni} the
following theorem was proved, see Proposition $22$.

\medskip

\t\label{basis of H^1(T)^s''}{\it Assume that
$H^1({\mathbb{CP}}^1,{\mathcal T}_2)^{\mathfrak s''}\ne 0$. Then we
have the following possibilities:

 \noindent $1)\;(k_1,k_2,k_3)=(2,2,2)$ and the basis of
$H^1({\mathbb{CP}}^1,{\mathcal T}_2)^{\mathfrak s''}$ is given by
\begin{equation}\label{basis s''1}
[\frac{1}{x^3}\xi_1\xi_2\xi_3\frac{\partial}{\partial \xi_3}-
\frac{1}{2x^2}\xi_1\xi_2\xi_3\frac{\partial}{\partial \xi_2}+
\frac{1}{2x}\xi_1\xi_2\xi_3\frac{\partial}{\partial \xi_1}];
\end{equation}

\noindent  $2)\;(k_1,k_2,k_3)=(3,3,3)$ and the basis of
$H^1({\mathbb{CP}}^1,{\mathcal T}_2)^{\mathfrak s''}$ is given by
\begin{equation}\label{basis s''2}
\begin{array}{c}
[\frac{1}{x^3}\xi_1\xi_2\frac{\partial}{\partial x}+
\frac{1}{2x^2}\xi_1\xi_3\frac{\partial}{\partial x}+
\frac{1}{2x}\xi_2\xi_3\frac{\partial}{\partial x}\\

+\frac{3}{8x^2}\xi_1\xi_2\xi_3\frac{\partial}{\partial \xi_1}-
\frac{3}{4x^3}\xi_1\xi_2\xi_3\frac{\partial}{\partial \xi_2}+
\frac{9}{4x^4}\xi_1\xi_2\xi_3\frac{\partial}{\partial \xi_3}].
\end{array}
\end{equation}
}

\bigskip

{\bf 5. Classification of even-homogeneous supermanifolds}

In Section $4$ we calculated a basis of
$H^1({\mathbb{CP}}^1,{\mathcal T}_2)^{\mathfrak s}$ and
$H^1({\mathbb{CP}}^1,{\mathcal T}_2)^{\mathfrak s'}$ and gave a
basis of $H^1({\mathbb{CP}}^1,{\mathcal T}_2)^{\mathfrak s''}$,
which were calculated in \cite{Bunegina_Oni}. In this section we
will complete the classification of even-homogeneous supermanifolds,
i.e. we will find out, which vectors of these spaces belong to
different orbits of the action of
 $\operatorname{Aut}\bold E$ on
$H^1({\mathbb{CP}}^1,{\mathcal T}_2)$.

Let $(\xi_i)$ be a local basis of $\bold E$ over $U_0$ and $A$ be an
automorphism of $\bold E$. Assume that $A(\xi_j)=\sum
a_{ij}(x)\xi_i$. In $U_1$ we have
$$
A(\eta_j) = A(y^{k_j}\xi_j) = \sum y^{k_j- k_i}a_{ij}(y^{-1})\eta_i.
$$
Therefore, $a_{ij}(x)$ is a polynomial in $x$ of degree no greater
than $k_j-k_i$, if
 $k_j-k_i\ge 0$ and $0$, if $k_j-k_i<0$.
We will denote by $b_{ij}$ the entries of the matrix $B=A^{-1}$. The
entries are also polynomials in $x$ of degree no greater than
$k_j-k_i$. We will need the following formulas:
\begin{equation}\label{A(vector field)}
\begin{array}{c}
A(\xi_1\xi_2\xi_3\frac{\partial}{\partial \xi_k})A^{-1} =
\det(A)\sum_s
b_{ks} \xi_1\xi_2\xi_3 \frac{\partial}{\partial \xi_s};\\
\rule{0pt}{6mm}A(\xi_i\xi_j\frac{\partial}{\partial x})A^{-1} =
\det(A) \sum_{k<s} (-1)^{l+r}b_{lr}
\xi_k\xi_s\frac{\partial}{\partial x} + \\
+\det(A)\sum_s b'_{ls}\xi_i\xi_j\xi_l\frac{\partial}{\partial
\xi_s}.
\end{array}
\end{equation}
where $i < j$, $l\ne i,j$, $r\ne k,s$ and $b'_{ls} =
\frac{\partial}{\partial x}(b_{ls})$.

\medskip

\t\label{Theor class 1}[Classification of
$\mathfrak{s}$-even-homogeneous supermanifolds.]{\it

$1.$ Assume that
$$
\{ k_1,4-k_1,k_3\}\neq \{-2,0,4\}, \quad \{ k,2-k,0\}\neq \{-2,0,4\}
$$
as sets. Then there exists a unique up to isomorphism ${\mathfrak
s}$-even-homogeneous non-split supermanifold with retract

$\textbf{a.}$ $(k_1,4-k_1,k_3)$, which correspond to the cocycle
$$
\frac1x\xi_1\xi_2\frac{\partial}{\partial x}+
\frac{k_3}{2x^2}\xi_1\xi_2\xi_3\frac{\partial}{\partial \xi_3};
$$

$\textbf{b.}$ $(k,2-k,0)$, which correspond to the cocycle

$$
b)\; \frac{1}{x}\xi_1\xi_2\xi_3\frac{\partial}{\partial \xi_3}.
$$

$2.$ There exist two up to isomorphism ${\mathfrak
s}$-even-homogeneous non-split supermanifolds  with retract
$(-2,0,4)$. The corresponding cocycles are
$$
a)\; z=\frac{1}{x}\xi_1\xi_2\xi_3\frac{\partial}{\partial \xi_2},
\quad b)\; z=\frac{1}x\xi_2\xi_3\frac{\partial}{\partial x}-
\frac{1}{x^2}\xi_1\xi_2\xi_3\frac{\partial}{\partial \xi_1}.
$$ }

\medskip

\noindent{\it Proof.} Since $m=3$, the number of different pairs
$i<j$ is less than or equal to $3$. It follows from the Theorem
\ref{basis of H^1(T)^a} that $\dim H^1({\mathbb{CP}}^1,{\mathcal
T}_2)^{\mathfrak s}\le 3$. It is easy to see that $\dim
H^1({\mathbb{CP}}^1,{\mathcal T}_2)^{\mathfrak s}=3$ if and only if
$k_1=k_2=k_3=2$. Let us take $A\in \operatorname{Aut}\bold
E=\operatorname{GL}_3({\mathbb{C}})$. Recall that
$\operatorname{Int} A(z)=AzA^{-1}$. The direct calculation shows,
see (\ref{A(vector field)}), that in the basis
$$
\begin{array}{c}
v_1=\frac1x\xi_2\xi_3\frac{\partial}{\partial x}+
\frac1{x^2}\xi_1\xi_2\xi_3\frac{\partial}{\partial \xi_1}, \quad
v_2=-\frac1x\xi_1\xi_3\frac{\partial}{\partial x}+
\frac1{x^2}\xi_1\xi_2\xi_3\frac{\partial}{\partial \xi_2}\\
v_3=\frac1x\xi_1\xi_2\frac{\partial}{\partial x}+
\frac1{x^2}\xi_1\xi_2\xi_3\frac{\partial}{\partial \xi_3},
\end{array}
$$
 the
automorphism ${\operatorname{Int}}\, A$ is given by
\begin{equation}\label{Int A(v)=}
{\operatorname{Int}}\, A(v_i) = \det A\sum_jb_{ij}v_j.
\end{equation}
 Note that for any matrix $C\in
\operatorname{GL}_3(\Bbb C)$ there exists a matrix $B$ such that
$C=\frac1{\det\,B}B$. Indeed, we can put
$B=\frac{1}{\sqrt{\det\,C}}C$. Let us take a cocycle $z = \sum
\alpha_i v_i\in H^1({\mathbb{CP}}^1,{\mathcal T}_2)^{\mathfrak
s}\setminus \{0\}$. Obviously, it exists a matrix
$D\in\operatorname{GL}_3(\Bbb C)$ such that $D(z)=(0,0,1)$.
Therefore, in the case $(2,2,2)$ there exists a unique up to
isomorphism $\mathfrak{s}$-even-homogeneous non-split supermanifold
given by the cocycle $\frac1x\xi_1\xi_2\frac{\partial}{\partial x}+
\frac1{x^2}\xi_1\xi_2\xi_3\frac{\partial}{\partial \xi_3}$.


Assume now that  $\dim H^1({\mathbb{CP}}^1,{\mathcal
T}_2)^{\mathfrak s}=2$. Let us consider three cases.

\smallskip

\noindent \textbf{1.} Assume that $H^1({\mathbb{CP}}^1,{\mathcal
T}_2)^{\mathfrak s}$ is generated by two cocycles from the item $1$
of Theorem \ref{basis of H^1(T)^a}. Obviously, we may consider only
the case $k_1+k_2=4,\, k_1+k_3=4$. It follows that $k_2=k_3$. Denote
$k_2:=k\ne 2$. Let us take $z\in H^1({\mathbb{CP}}^1,{\mathcal
T}_2)^{\mathfrak s}\setminus \{0\}$. Then
$z=\frac{\alpha}{x}\xi_1\xi_2\frac{\partial}{\partial x}+
\frac{k\alpha}{2x^2}\xi_1\xi_2\xi_3\frac{\partial}{\partial \xi_3}+
\frac{\beta}{x}\xi_1\xi_3\frac{\partial}{\partial x}+
\frac{k\beta}{2x^2}\xi_1\xi_3\xi_2\frac{\partial}{\partial \xi_2}$.
The group $\operatorname{Aut}\bold E$ contains in this case the
subgroup $H$:
\begin{equation}\label{subgroup in Aut E1}
H:= \left\{ \left(
\begin{array}{ccc}
a_{11}&0&0\\
0&a_{22}&a_{23}\\
0&a_{32}&a_{33}\\
\end{array}
\right)\right\}.
\end{equation}
Let us take $A\in H$, denote
$v_1:=-\frac1x\xi_1\xi_3\frac{\partial}{\partial x}+
\frac{k}{2x^2}\xi_1\xi_2\xi_3\frac{\partial}{\partial \xi_2}$,
$v_2:=\frac1x\xi_1\xi_2\frac{\partial}{\partial x}+
\frac{k}{2x^2}\xi_1\xi_2\xi_3\frac{\partial}{\partial \xi_3}$. Using
(\ref{A(vector field)}) or (\ref{Int A(v)=}) we see that the
operator ${\operatorname{Int}}\, A$ is given in the basis $v_1$,
$v_2$ by:
$$
\det A \left(
\begin{array}{cc}
b_{22}&b_{23}\\
b_{32}&b_{33}\\
\end{array}
\right).
$$
Obviously, for any cocycle $z=(-\beta,\alpha)\ne 0$ there exists a
matrix $C\in \GL_3(\mathbb{C})$ such that
 $C(z)=(0,1)$. Therefore, in the case $(4-k,k,k)$ there exists a
unique up to isomorphism $\mathfrak{s}$-even-homogeneous non-split
supermanifold given by the cocycle
$\frac1x\xi_1\xi_2\frac{\partial}{\partial x}+
\frac{k}{2x^2}\xi_1\xi_2\xi_3\frac{\partial}{\partial \xi_3}$.

\smallskip

\noindent \textbf{2.} Assume that $H^1({\mathbb{CP}}^1,{\mathcal
T}_2)^{\mathfrak s}$ is generated by two cocycles from the item $2$
of Theorem \ref{basis of H^1(T)^a}. We may consider only the case $
k_1+k_2=2,\,k_1+k_3=2,\, k_2=k_3=0$. It follows that $(k_1,k_2,k_3)
= (2,0,0)$. Let us take $z\in H^1({\mathbb{CP}}^1,{\mathcal
T}_2)^{\mathfrak s}\setminus \{0\}$. Then
$z=\frac{\alpha}x\xi_1\xi_2\xi_3\frac{\partial}{\partial \xi_3}+
\frac{\beta}{x}\xi_1\xi_2\xi_3\frac{\partial}{\partial \xi_2}$,
where $\alpha,\beta\in \mathbb{C}$. As above, the group
$\operatorname{Aut}\bold E$ contains the subgroup $H$ given by
$(\ref{subgroup in Aut E1})$. As above using the basis
$v_1=\frac1{x}\xi_1\xi_2\xi_3\frac{\partial}{\partial \xi_2}$,
$v_2=\frac1{x}\xi_1\xi_2\xi_3\frac{\partial}{\partial \xi_3}$, we
show that in the case $(2,0,0)$ there exists a unique up to
isomorphism $\mathfrak{s}$-even-homogeneous non-split supermanifold
given by the cocycle
$\frac1{x}\xi_1\xi_2\xi_3\frac{\partial}{\partial \xi_2}$.


\smallskip

\noindent \textbf{3.} Assume that $H^1({\mathbb{CP}}^1,{\mathcal
T}_2)^{\mathfrak s}$ is generated by one cocycle from the item $1$
and by one cocycle from the item $2$ of Theorem $4$. We may consider
only the case $k_2+k_3=4,\, k_1+k_3=2,\,k_2=0$, i.e. $(k_1,k_2,k_3)
= (-2,0,4)$. Let us take $z\in H^1({\mathbb{CP}}^1,{\mathcal
T}_2)^{\mathfrak s}\setminus \{0\}$. Then
$\,z=\frac{\alpha}x\xi_2\xi_3\frac{\partial}{\partial x}-
\frac{\alpha}{x^2}\xi_1\xi_2\xi_3\frac{\partial}{\partial \xi_1}+
\frac{\beta}{x}\xi_1\xi_2\xi_3\frac{\partial}{\partial \xi_2}$ for
certain $\alpha,\beta\in \mathbb{C}$. Let us take $A\in
\operatorname{Aut} \bold E$. Using Theorem \ref{basis of H^1} and
\ref{A(vector field)}, we get
$$
\begin{array}{rcl}
A([\frac1x\xi_2\xi_3\frac{\partial}{\partial
x}])A^{-1}&=&[b_{11}\det A(\frac1x\xi_2\xi_3\frac{\partial}{\partial
x}
+(b_{12})'\frac1x\xi_1\xi_2\xi_3\frac{\partial}{\partial \xi_2})];\\
A([\frac1{x^2}\xi_1\xi_2\xi_3\frac{\partial}{\partial
\xi_1}])A^{-1}&=& [b_{11} \det
A\frac1{x^2}\xi_1\xi_2\xi_3\frac{\partial}{\partial \xi_1}];\\
A([\frac1{x}\xi_1\xi_2\xi_3\frac{\partial}{\partial
\xi_2}])A^{-1}&=& [\det A
(b_{22}\frac1{x}\xi_1\xi_2\xi_3\frac{\partial}{\partial \xi_2})],
\end{array}
$$
where $(b_{12})':= \frac{\partial}{\partial x}(b_{12})$. Consider
the subgroup $H\!=\{\!\operatorname{diag}(a_{11},a_{22},a_{33})\}$
of $\operatorname{Aut} \bold E$. Let us choose the basis
$v_1=\frac1x\xi_2\xi_3\frac{\partial}{\partial x}-
\frac1{x^2}\xi_1\xi_2\xi_3\frac{\partial}{\partial \xi_1}$,
$v_2=\frac1{x}\xi_1\xi_2\xi_3\frac{\partial}{\partial \xi_2}$ and
take $A\in H$. Then the operator ${\operatorname{Int}}\,A$ is given
by the matrix
$$
(\det A) \operatorname{diag}(b_{11},b_{22})
$$ in the
basis $v_1$, $v_2$. Obviously, for any cocycle $z=(\alpha,\beta)\ne
0$ there exists an operator ${\operatorname{Int}}\,A$ such that:
${\operatorname{Int}}\,A(z)=(1,1)$, if $\alpha\ne 0$, $\beta\ne 0$,
${\operatorname{Int}}\,A(z)=(0,1)$, if $\alpha=0$, $\beta\ne 0$,
${\operatorname{Int}}\,A(z)=(1,0)$, if $\alpha\ne 0$, $\beta=0$. Let
us take
$$
A= \left(
     \begin{array}{ccc}
       1 & -x & 0 \\
       0& 1 & 0 \\
       0 & 0 & 1 \\
     \end{array}
   \right)\in \operatorname{Aut} \bold E.
$$
The direct calculation shows that $A(v_1)A^{-1} = v_1+v_2$. In other
words, $v_1$ and $v_1+v_2$ corresponds to one orbit of the action
${\operatorname{Int}}$. Since $b_{11}\ne 0$, we see that the
cocycles  $(0,1)$ and $(1,0)$ correspond to different orbits of the
action ${\operatorname{Int}}$.

In the case $\dim H^1({\mathbb{CP}}^1,{\mathcal T}_2)^{\mathfrak
s}=1$ we may use the following proposition proved in \cite{OniNA}.

\medskip

\l\label{Prop a and ca}{\it If $\gamma \in
H^1({\mathbb{CP}}^1,{\mathcal T}_2)$ и $c\in {\mathbb{C}}\setminus
\{0\}$, then $\gamma$ and $c\gamma$ correspond to isomorphic
supermanifolds.}

\medskip

Theorem \ref{Theor class 1} follows.$\Box$

\medskip

\t\label{Theor class 2}[Classification of
$\mathfrak{s}'$-even-homogeneous supermanifolds.] {\it $1.$ There
exist two up to isomorphism ${\mathfrak s}'$-even-homogeneous
non-split supermanifolds with retract

\noindent  $a)$ $(2,2,1)$, which correspond to the cocycles
$$
\frac1x\xi_1\xi_2\frac{\partial}{\partial x}+
\frac{1}{2x^2}\xi_1\xi_2\xi_3\frac{\partial}{\partial \xi_3},\quad
\frac{1}{x^2}\xi_1\xi_2\xi_3\frac{\partial}{\partial
\xi_2}-\frac{1}{x}\xi_1\xi_2\xi_3\frac{\partial}{\partial \xi_1},
$$

\noindent $b)\;(2,2,3)$, which correspond to the cocycles
$$
\begin{array}{c}
\frac1x\xi_1\xi_2\frac{\partial}{\partial x}+
\frac{3}{2x^2}\xi_1\xi_2\xi_3\frac{\partial}{\partial \xi_3},\\
\frac{1}{x}\xi_2\xi_3\frac{\partial}{\partial x}+
\frac{1}{x^2}\xi_1\xi_3\frac{\partial}{\partial x}+
\frac{2}{3x^2}\xi_1\xi_2\xi_3\frac{\partial}{\partial \xi_1}
-\frac{4}{3x^3}\xi_1\xi_2\xi_3\frac{\partial}{\partial \xi_2}.
\end{array}
$$

$2.$ There exists a unique up to isomorphism ${\mathfrak
s}'$-even-homogeneous non-split supermanifold with retract

\noindent $a)\;(2,2,k)$, $k\ne 1,3$, which corresponds to the
cocycle
$$
\frac1x\xi_1\xi_2\frac{\partial}{\partial x}+
\frac{k}{2x^2}\xi_1\xi_2\xi_3\frac{\partial}{\partial \xi_3},
$$

\noindent $b)$ $(k,k,3-k)$, $k\ne 2$, which corresponds to the
cocycle
$$
\frac{1}{x^2}\xi_1\xi_2\xi_3\frac{\partial}{\partial \xi_2}-
\frac{1}{x}\xi_1\xi_2\xi_3\frac{\partial}{\partial \xi_1},
$$

\noindent $c)\; (k,k,5-k)$, $k\ne 2$, which corresponds to the
cocycle
$$
\frac{1}{x}\xi_2\xi_3\frac{\partial}{\partial x}+
\frac{1}{x^2}\xi_1\xi_3\frac{\partial}{\partial x}+
\frac{k}{3x^2}\xi_1\xi_2\xi_3\frac{\partial}{\partial \xi_1}
-\frac{2k}{3x^3}\xi_1\xi_2\xi_3\frac{\partial}{\partial \xi_2}.
$$

\noindent $d)\; (1,1,0)$, which corresponds to the cocycle
$$
\frac{1}{x}\xi_1\xi_2\xi_3\frac{\partial}{\partial \xi_3}.
$$

}

\medskip
\noindent{\it Proof.} By Theorem \ref{basis of H^1(T)^s'} and
Proposition \ref{Prop a and ca} we get $2.$

Let us prove $1.a$ Denote by $z$ a linear combination of cocycles
(\ref{basis s'1}). Let us take $A\in \operatorname{Aut} \bold E$.
Using (\ref{A(vector field)}), we get:
$$
\begin{array}{l}
A([\frac1x\xi_1\xi_2\frac{\partial}{\partial x}])A^{-1}= [\det
A(b_{33}\frac1x\xi_1\xi_2\frac{\partial}{\partial x} +
(b_{31})'\frac1x\xi_1\xi_2\xi_3\frac{\partial}{\partial \xi_1}\!+\!
(b_{32})'\frac1x\xi_1\xi_2\xi_3\frac{\partial}{\partial
\xi_2})];\\
\end{array}
$$
$$
\begin{array}{l}
A([\frac1{x^2}\xi_1\xi_2\xi_3\frac{\partial}{\partial
\xi_3}])A^{-1}\!=\! [\det
A(b_{33}\frac1{x^2}\xi_1\xi_2\xi_3\frac{\partial}{\partial
\xi_3}\!+\!
b_{32}\frac1{x^2}\xi_1\xi_2\xi_3\frac{\partial}{\partial \xi_2}+\\
+b_{31}\frac1{x^2}\xi_1\xi_2\xi_3\frac{\partial}{\partial \xi_1})];\\
A([\frac1{x^2}\xi_1\xi_2\xi_3\frac{\partial}{\partial
\xi_2}])A^{-1}\!\!=\! [\det
A(b_{21}\frac1{x^2}\xi_1\xi_2\xi_3\frac{\partial}{\partial
\xi_1}\!+\!
b_{22}\frac1{x^2}\xi_1\xi_2\xi_3\frac{\partial}{\partial \xi_2})],\\
A([\frac1{x}\xi_1\xi_2\xi_3\frac{\partial}{\partial
\xi_1}])A^{-1}\!\!=\! [\det
A(b_{11}\frac1{x}\xi_1\xi_2\xi_3\frac{\partial}{\partial \xi_1}\!+\!
b_{12}\frac1{x}\xi_1\xi_2\xi_3\frac{\partial}{\partial \xi_2})].
\end{array}
$$
Consider the subgroup
$H\!=\!\{\operatorname{diag}(a_{11},a_{11},a_{33},)\}$ of
$\operatorname{Aut} \bold E$ and $A\in H$. Again a direct
calculation shows that in the basis
$v_1\!=\!\frac1x\xi_1\xi_2\frac{\partial}{\partial x}+
\frac1{2x^2}\xi_1\xi_2\xi_3\frac{\partial}{\partial \xi_3}$,
$v_2\!=\!\frac1{x^2}\xi_1\xi_2\xi_3\frac{\partial}{\partial \xi_2}-
\frac1{x}\xi_1\xi_2\xi_3\frac{\partial}{\partial \xi_1}$ the
automorphism ${\operatorname{Int}}\,A$ is given by $(\det A)
\operatorname{diag}(b_{33},b_{11})$. Clearly, for
$z=(\alpha,\beta)\ne 0$, there exist an operator
${\operatorname{Int}}\,A$ such that:
${\operatorname{Int}}\,A(z)=(1,1)$, if $\alpha\ne 0$, $\beta\ne 0$,
${\operatorname{Int}}\,A(z)=(0,1)$, if $\alpha=0$, $\beta\ne 0$,
${\operatorname{Int}}\,A(z)=(1,0)$, if $\alpha\ne 0$, $\beta=0$.

Let us take
$$
A=\left(
\begin{array}{ccc}
1&0&0\\
0&1&1\\
-\frac23 x& 2 &1
\end{array} \right),
$$
A direct calculation shows that $A(v_1+v_2)A^{-1} = v_1$. Since
$b_{33}\ne 0$, we see that the cocycles  $(0,1)$ and $(1,0)$
correspond to different orbits of the action ${\operatorname{Int}}$.
We have got $1a)$. The proof of $1b)$ is similar. The result
follows.$\Box$

\medskip

\t\label{Theor class 3}[Classification of
$\mathfrak{s}''$-even-homogeneous supermanifolds.] {\sl There exist
a unique up to isomorphism ${\mathfrak s}''$-even-homogeneous
non-split supermanifold with retract $(2,2,2)$, which corresponds to
the cocycle $(\ref{basis s''1})$; and with retract $(3,3,3)$, which
corresponds to the cocycle $(\ref{basis s''2})$. }

\medskip
\noindent{\it Proof.} It follows from Theorem \ref{basis of
H^1(T)^s''} and Proposition \ref{Prop a and ca}.$\Box$

\medskip

Comparing Theorems \ref{Theor class 1}, \ref{Theor class 2} and
\ref{Theor class 3}, we get our main result:

\medskip

\t\label{Theor class obshch}[Classification of even-homogeneous
supermanifolds.] {\sl

$1.$ There exist two up to isomorphism even-homogeneous non-split
supermanifolds with retract

\noindent  $a)$ $(2,2,1)$, which correspond to the cocycles
$$
\frac1x\xi_1\xi_2\frac{\partial}{\partial x}+
\frac{1}{2x^2}\xi_1\xi_2\xi_3\frac{\partial}{\partial \xi_3},\quad
\frac{1}{x^2}\xi_1\xi_2\xi_3\frac{\partial}{\partial
\xi_2}-\frac{1}{x}\xi_1\xi_2\xi_3\frac{\partial}{\partial \xi_1};
$$

\noindent $b)\;(2,2,3)$, which correspond to the cocycles
$$
\begin{array}{c}
\frac1x\xi_1\xi_2\frac{\partial}{\partial x}+
\frac{3}{2x^2}\xi_1\xi_2\xi_3\frac{\partial}{\partial \xi_3},\\
\frac{1}{x}\xi_2\xi_3\frac{\partial}{\partial x}+
\frac{1}{x^2}\xi_1\xi_3\frac{\partial}{\partial x}+
\frac{2}{3x^2}\xi_1\xi_2\xi_3\frac{\partial}{\partial \xi_1}
-\frac{4}{3x^3}\xi_1\xi_2\xi_3\frac{\partial}{\partial \xi_2};
\end{array}
$$

\noindent $c)\;(2,2,2)$, which correspond to the cocycles
$$
\begin{array}{c}
\frac1x\xi_1\xi_2\frac{\partial}{\partial x}+
\frac{1}{2x^2}\xi_1\xi_2\xi_3\frac{\partial}{\partial \xi_3},\\
\frac{1}{x^3}\xi_1\xi_2\xi_3\frac{\partial}{\partial \xi_3}-
\frac{1}{2x^2}\xi_1\xi_2\xi_3\frac{\partial}{\partial \xi_2}+
\frac{1}{2x}\xi_1\xi_2\xi_3\frac{\partial}{\partial \xi_1};
\end{array}
$$

\noindent $d)\;(-2,0,4)$, which correspond to the cocycles
$$
\frac{1}{x}\xi_1\xi_2\xi_3\frac{\partial}{\partial \xi_2}, \quad
\frac{1}x\xi_2\xi_3\frac{\partial}{\partial x}-
\frac{1}{x^2}\xi_1\xi_2\xi_3\frac{\partial}{\partial \xi_1}.
$$

\noindent $2.$ $a)$ Assume that
$$
\{k,4-k,k_3\}\ne \{-2,0,4\}, \, \{2,2,1\},\, \{2,2,3\},\, \{2,2,2\}.
$$
Then there exists a unique up to isomorphism even-homogeneous
non-split supermanifold corresponding to the cocycle
$$
\frac1x\xi_1\xi_2\frac{\partial}{\partial x}+
\frac{k_3}{2x^2}\xi_1\xi_2\xi_3\frac{\partial}{\partial \xi_3}.
$$

\noindent $b)$ Assume that
$$
\{k,2-k,0\}\ne \{-2,0,4\}.
$$
Then there exists a unique up to isomorphism even-homogeneous
non-split supermanifold corresponding to the cocycle
$$
\frac{1}{x}\xi_1\xi_2\xi_3\frac{\partial}{\partial \xi_3}.
$$

 There exists a unique up to isomorphism even-homogeneous non-split supermanifold with retract

\noindent $c)$ $(k,k,3-k)$, $k\ne 2$, which corresponds to the
cocycle
$$
\frac{1}{x^2}\xi_1\xi_2\xi_3\frac{\partial}{\partial \xi_2}-
\frac{1}{x}\xi_1\xi_2\xi_3\frac{\partial}{\partial \xi_1},
$$

\noindent $d)\; (k,k,5-k)$, $k\ne 2$, which corresponds to the
cocycle
$$
\frac{1}{x}\xi_2\xi_3\frac{\partial}{\partial x}+
\frac{1}{x^2}\xi_1\xi_3\frac{\partial}{\partial x}+
\frac{k}{3x^2}\xi_1\xi_2\xi_3\frac{\partial}{\partial \xi_1}
-\frac{2k}{3x^3}\xi_1\xi_2\xi_3\frac{\partial}{\partial \xi_2}.
$$

\noindent $d)\; (3,3,3)$, which corresponds to the cocycle
(\ref{basis s''2}). $\Box$

 }

\medskip

By the similar argument as in \cite{Bunegina_Oni}, Corollary of
Theorem $1$, we get:

\noindent{\bf Corollary.} Any non-split even-homogeneous
supermanifold $\mathcal{M}$ of dimension $1|2$, where $\mathcal{M}_0
= \mathbb{CP}^1$, is isomorphic to $\mathbb{Q}^{1|2}$.
\medskip

Here $\mathbb{Q}^{1|2}$ is the (homogeneous) supermanifold
corresponding to the cocycle
$x^{-1}\xi_1\xi_2\frac{\partial}{\partial x}$ (see
\cite{Bunegina_Oni} for more details).

\medskip

\noindent{\bf Remark 1.} Theorem \ref{Theor class obshch}
 gives rise to a description of even-homogeneous supermanifolds in
 terms of local charts and coordinates. Indeed, let $\mathcal{M}$ be
 any supermanifold of dimension $1|m$, $m\leq 3$, with underlying space $\mathbb{CP}^1$,
 $v$ be the corresponding cocycle by Theorem \ref{Green corollary} and $(U_0,
 \mathcal{O}_{\gr\mathcal{M}}|_{U_0})$, $(U_1,
 \mathcal{O}_{\gr\mathcal{M}}|_{U_1})$ be two standard charts of the retract $\gr\mathcal{M}$ with
 coordinates $(x, \xi_1,\xi_2, \xi_3)$ and $(y, \eta_1,\eta_2,
 \eta_3)$, respectively. In $U_0\cap U_1$ we have:
 $$
 y=x^{-1},\quad \eta_i = x^{-k_i}\xi_i, \,\,i=1,2,3.
 $$

Consider an atlas on $\mathcal{M}$: $(U_0,
 \mathcal{O}_{\mathcal{M}}|_{U_0})$, $(U_1,
 \mathcal{O}_{\mathcal{M}}|_{U_1})$, with coordinates $(x', \xi'_1,\xi'_2, \xi'_3)$ and $(y', \eta'_1,\eta'_2,
 \eta'_3)$, respectively. Then the transition  function
of $\mathcal{M}$ in $U_0\cap U_1$ have the form
$$
y'=(\id +v)(x'^{-1}),\quad \eta_i = (\id +v)((x')^{-k_i}\xi'_i),
\,\,i=1,2,3.
$$

\medskip

\noindent{\bf Remark 2.} The supermanifold $\mathcal{M}$ with the
retract $(k,2-k,0)$, corresponding to the cocycle
$\frac{1}{x}\xi_1\xi_2\xi_3\frac{\partial}{\partial \xi_3}$, which
was lost in \cite{Bunegina_Oni}, in even-homogeneous but not
homogeneous. Hence the main result in \cite{Bunegina_Oni}, Theorem
$1$,  is correct.


\noindent{\it Elizaveta Vishnyakova}

\noindent {Mathematisches Forschungsinstitut Oberwolfach and}

\noindent University of Luxembourg

 \noindent {\emph{E-mail address:}
\verb"VishnyakovaE@googlemail.com"}

\end{document}